\documentclass[11pt]{amsart}
\usepackage{amsaddr, amsmath,amsfonts,amssymb,amsthm,bbm}
\usepackage{mathrsfs}
\usepackage{float}
\usepackage{graphicx,color,cases,mathtools}
\usepackage{subfig}%
\newtheorem{Qu}{Question}
\def\beQ{\begin{Qu}} \def\eeQ{\end{Qu}}

 \def\ssec{\subsection}
 
 \def\wk{well-known}

 \def\procs{processes} 
\def\Rui{\Psi} \def\sRui{\ovl{\Rui}}
 \def\fp{first passage }
\def\ovl{\overline}
\def\lev{L\'evy }    
\def\dd{draw-down } \def\sn{spectrally negative }
\long\def\symbolfootnote[#1]#2{
\begingroup
\def\thefootnote{\fnsymbol{footnote}}\footnote[#1]{#2}
\endgroup} \def\I{\infty} \def\Eq{\Leftrightarrow}
\def\fn{\symbolfootnote} \def\und{\underline}  \def\T{\widetilde}
\def\CL{Cram\'er-Lundberg }  
   
\def\BEN{\begin{enumerate}}  \def\BI{\begin{itemize}}
\def\EEN{\end{enumerate}}   \def\EI{\end{itemize}} \def\im{\item} \def\Lra{\Longrightarrow} \def\D{\Delta} \def\eqr{\eqref}
\def\no{\nonumber} 

   \def\de{\delta}  \def\b{\beta}
  \def\th{\theta} 
\def\e{\epsilon} \newcommand{\ol}{\overline}
\def\k{\kappa} \def\l{\lambda} \def\a{\alpha} 
  
  \def\ith{it holds that }
  \def\r{\rho} \def\s{\sigma} \def\F{\Phi} \def\f{\varphi}\def\L{L} \def\U{U}
  \def\bc{\begin{cases}
  }    
\def\ec{\end{cases}} 
  \def\qu{\quad} \def\for{\forall}
  \def\beq{\begin{eqnarray}} \def\eeq{\end{eqnarray}}
   \def\be{\begin{equation}} \def\ee{\end{equation}}
\def\bea{\begin{eqnarray*}}
\def\eea{\end{eqnarray*}} \def\la{\label}
  
 \def\q{q} 
\def\le{\left} \def\ri{\right}
 \def\fr{\frac} \def\Y{Y}  \def\t{\tau} \def\ta{T_{a,-}}  \def\tb{T_{b,+}}
\def\tz{T_{0,-}} \def\sec{\section} \def \fe{for example } 
 \def\deF{de Finetti } 

\newtheorem{Thm}{Theorem}
\newtheorem{Ass}{Assumption}
\newtheorem{Lem}{Lemma}
\newtheorem{Exa}[Lem]{Example}

\newtheorem{Cor}[Lem]{Corollary}
\newtheorem{Def}[Lem]{Definition}
\newtheorem{Rem}[Lem]{Remark}
\newtheorem{Exe}[Lem]{Exercice}
\newtheorem{Con}{Conjecture}
\def\beXe{\begin{Exe}} \def\eeXe{\end{Exe}}
\def\eeD{\end{Def}} \def\beD{\begin{Def}}
\def\beXa{\begin{Exa}} \def\eeXa{\end{Exa}}
\def\beR{\begin{Rem}} \def\eeR{\end{Rem}}
\def\beL{\begin{Lem}} \def\eeL{\end{Lem}}
\newtheorem{Pro}[Lem]{Proposition}
\def\beP{\begin{Pro}} \def\eeP{\end{Pro}}
\def\beC{\begin{Cor}} 
\def\beT{\begin{Thm}}
  \def\eeT{\end{Thm}}   \def\Mar{Markov }
\def\eeC{\end{Cor}}   \def\sur{survival }
 \def\cP{compound Poisson }
 
  \def\gen{generalized }

\newcommand{\pd}[2]{\frac{\partial #1}{\partial #2}}

\def\prob{problem} \def\probs{problems}  \def\str{strong Markov property}
\def\frt{furthermore }  \def\td{\tau}
\def\strs{strong Markov processes}    \def\AY{Azema-Yor }
 \def\OU{Ornstein-Uhlenbeck } 
  \def\dif{diffusion}
\def\1{\mathbf{1}}  \newcommand{\md}{\mathrm d}

     \def\Thr{Therefore, }
    \def\iffa{} 
    
    \def\sats{satisfies}  \def\fun{function }  \def\thr{therefore } \def\wf{we find that } \def\inc{increasing}
   
\providecommand{\norm}[1]{\left\lVert#1\right\rVert}
\providecommand{\abs}[1]{\left\lvert#1\right\rvert}
\providecommand{\pr}[1]{\left(#1\right)} 
\providecommand{\pp}[1]{\left[#1\right]} 
\providecommand{\set}[1]{\left\lbrace#1\right\rbrace} 
\providecommand{\scal}[1]{\left\langle#1\right\rangle}

\newcommand{\figu}[3]{
\begin{figure}[!h]
\begin{center}
{\includegraphics[width=13 cm, height=8 cm]{#1}}
\end{center}

\vspace{-0.2cm}
\caption{\hspace{0.25cm}#2\label{f:#1}}
\end{figure}
}
\begin{document}
\title[Optimal Draw-Down]{A Pontryaghin maximum principle approach  for the optimization of dividends/consumption of spectrally negative Markov
processes, until a generalized draw-down time}

\author{Florin Avram}  \address{Laboratoire de Math\'ematiques Appliqu\'ees, Universit\'e de Pau,  France} \email{ florin.avram@univ-Pau.fr}
\author{Dan Goreac}\address{Universit\'e Paris-Est, LAMA (UMR 8050), UPEMLV, UPEC, CNRS, F-77454, Marne-la-Vall\'ee, France}\email{dan.goreac@u-pem.fr}

\date{\today}

\begin{abstract}

The first motivation of our  paper  is to explore further the idea
that, in risk control  problems, it may be profitable to base decisions both on the position of the underlying process $X_t$ and on its supremum $\ovl X_t:=\sup_{0\leq s\leq t} X_s$. Strongly connected to Azema-Yor/generalized draw-down/trailing stop time
(see \cite{AY}), this framework provides a natural unification
of  \dd and classic \fp times.

We illustrate here the potential of this unified framework by solving a variation of the De Finetti problem of maximizing expected discounted cumulative dividends/consumption gained under a barrier policy, until an optimally chosen Azema-Yor time, with a general spectrally negative Markov model.

While previously studied cases of this problem \cite{APP,SLG,alvarez1998optimal,AVZ,hening2018optimal,WangZhou} assumed  either  \lev  or  \dif \; models, and the \dd function to be fixed, we describe, for a general spectrally negative Markov model, not only the optimal barrier but also the optimal draw-down function. This is achieved by solving a variational problem tackled by Pontryaghin’s maximum principle. As a by-product we show that in the \lev case the classic first passage solution is indeed optimal; in the diffusion case, we obtain the optimality equations, but the  existence of solutions improving the classic ones is left for future work.

\end{abstract}
\maketitle


\tableofcontents

\sec{A brief review of the classic \sn \fp theory 
\la{s:r}}

{\bf Control of dividends/optimal consumption and capital injections}. Many  control problems in risk theory  concern versions  of  $X_t$ which are
 reflected/constrained/regulated at first passage times (below or above):
 \begin{eqnarray} \la{Xd}&& X_t^{[a}=X_t  + \; \L_t, \qu  X^{b]}_t=X_t   -  \U_t. \la{ref}\end{eqnarray}
Here, \beq \la{RR} \no && \L_t=\L^{[a }_t=-(\und X_t-a)_-, \quad \und X_t=\inf_{0 \leq t'\leq t} X_{t'},\\ \no && \U_t=\U^{b]}_t=\le(\ovl X_t -b\ri)_+, \quad \ovl X_t := \sup_{0 \leq t'\leq t} X_{t'},\eeq
  are  the minimal  'Skorohod regulators' constraining $X_t$ to be larger than $a$, and    smaller than $b$, respectively, and we use the notation $x_+=\max(x,0)$ and $x_-=\min(x,0)$.

 Financial management and other applications require also  studying
 the running maximum and  the process reflected at its maximum/drawdown $$Y_t=\ovl{ X}_t-X_t, \quad \ovl X_t= \sup_{0 \leq t'\leq t} X_s,$$
as well as the running infimum and  the process reflected from below/drawup $$\und Y_t= X_t-\und X_t, \quad \und X_t=\inf_{0 \leq t'\leq t} X_s.
$$

 The first passage times of the reflected processes, called
 draw-down/regret time and  draw-up time, respectively, are defined for $d>0$ by
\begin{equation} \label{dd}
\begin{aligned}
\t_{d}&:=\inf\{t\geq 0: \ovl X(t)- X(t) \geq d \},\\
\und \t_{d}&:=\inf\{t\geq 0:  X(t)- \und X(t) \geq d \}.
\end{aligned}
\end{equation}

 \iffa
Such times  
turn out to be optimal in several stopping problems, in  statistics \cite{page1954continuous}    in mathematical finance/risk theory   and in queueing. More specifically, they figure in risk theory  problems involving  dividends at a fixed barrier or capital injections, and in studying idle times until a buffer reaches capacity in queueing theory--  see for example \cite{taylor1975stopped,Leh,shepp1993russian,AKP,MP,shiryaevthou,Carr,LLZ17,LLZ17b}--
 for results and references to  numerous applications of draw-downs and draw-ups.

{\bf  Optimization of dividends}.
One  important optimization problem, going back  to  Bruno de Finetti \cite{deF}, is  to
  estimate the maximal expected discounted cumulative dividends of a financial company until its ruin.

   Solving this   is nontrivial even for  \lev process  {without positive jumps} (a ``multi-band" continuation region may be necessary); therefore,
given the usual data uncertainty inherent to real problems, it is  reasonable to restrict to simpler dividends  policies which distribute all surpluses above a fixed level $b$, called {\bf dividends barrier}.
For fixed $b$, we arrive then to
{\bf the  optimal dividends  barrier problem with ruin stopping}.

Note that the  time of ruin of the ``Skorohod regulated process"  $X^{b]}_t= X_t  -  \U_t$
    may be decomposed as:
\beq \la{Trb} T_0^{b]} :=\tz \; {\mathbf 1}_{\tz< \tb}+ \t_b \;  {\mathbf 1}_{\tb<\tz}.\eeq

The  ``\deF barrier problem" consists in maximizing over $b$
the present value of all dividend payments  at the barrier  $b$, until the  time  $T_0^{b]}$:
 \be \la{tdiv}
 V^{b}(x)=V_{\q}^{b}(x):=  E_x \left[ \int_{0}^{T_0^{b]} } e^{-q  t}d(\ol {X}_t-b)_+ \right].\ee

  {\bf  In the case of a \sn \lev process} $X_t$, the value function \eqr{tdiv} and  many other results may be expressed in terms of the
  scale function $W_q$ \cite{Ber,Kyp}.
In preparation for
 the \sn Markov case, we will also  express the value function in terms of the logarithmic derivative \beq \la{nu} \nu_{\q}(x):=\fr{W_{\q}'(x+)}{W_{\q}(x)}\eeq
 since it became apparent in \cite{LLZ17b,ALL} that for \sn
\Mar \procs \; it is more convenient  to introduce first
a natural extension
$\nu_\q(x,y)$  (defined via the limit  \eqr{e:b}) of the logarithmic derivative  than the corresponding extension $W_\q(x,y)$.

The
 \lev factorization of the ``gambler's winning/survival" probability may also be written as
 \beq \la{Wnu}   E_{x}\left[  e^{-q\tb
}1_{\left\{  \tb<\ta\right\}  }\right]  =\fr{W_{\q}(x-a)}{W_{\q}(b-a)}=e^{-\int_{x}^{b} \nu_q(s-a)ds}.\eeq

Applying now the \str \; in \eqr{tdiv} yields
 \bea
 V^{b}(x)=E_{x}\left[  e^{-q\tb} ;  \tb<\tz   \right] E_b \left[ \int_{0}^{T_0^{b]} } e^{-q  t}d(\ol {X}_t-b) \right]\eea
 \be = \fr{W_{\q}(x)}{W_{\q}(b)} V^{b}(b)=\fr{W_{\q}(x)}{W_{\q}(b)} \fr{W_{\q}(b)}{W_{\q}'(b)}=e^{-\int_{x}^{b} \nu_\q(m)dm}\fr 1{\nu_\q(b)},\la{VLev}\ee
 where   we have used \eqr{Wnu}  and
 $$V^{b}(b)=\fr 1{\nu_\q(b)}$$
 cf. \cite{Ber,kyprianou2007distributional,albrecher2018linking}.
 To understand the last equality, note that the  dividends starting from $b$ equal the local time spent at the reflecting boundary $b$, and that the latter has an  exponential law, the rate  of which is $\nu_{\q}(b)$, since the  function $\nu_{\q}(x)$ is the {rate of downwards excursions strictly larger than $x$, and occurring  before an exponential horizon of rate $q$} \cite{Ber,Doney}.

 We will make use below of the fact that $\nu_{\q}(x)$ is nonincreasing  and that
  \beq \la{nub} \nu_{\q}(x) \geq \F_\q \eeq where $ \F_\q$ is the unique positive root of the \CL equation
   \cite{Ber,Kyp}
   \beq \k(s):=Log \Big(E_0[e^{s X_1}]\Big)=q. \la{sym}\eeq

 \begin{Ass} \la{As} To be able to write below equations like $\nu_\q(x) =\frac{W_{\q}^ {\prime}(x)}{W_\q(  x)}$
  and  formulas like \eqr{costate}
   we will assume throughout the paper that
 $W_{\q}(x)$ is  three times
  differentiable in the \lev case. In the \sn Markov case, we will assume the scale function $W_{\q} (x,y)$ (see last section) to be
    three times
  differentiable in $x$, or,  alternatively,   $\nu_\q(x,y) :=\frac{\pd{W_{\q} (x,y)}{x}}{W_\q(  x,y)}$ will  be assumed to be twice differentiable.
\end{Ass}

 See \cite{chan2011smoothness} for more information on the smoothness of scale functions for \lev processes, and note this problem has not yet been studied for \sn \Mar \procs.

In conclusion, the \lev De-Finetti barrier objective
   has a simple expression  in terms of either the $W_\q $ scale function or of $\nu_{\q}$:
   \begin{equation} V^{b}(x)= \bc \frac{W_\q(  x)}{W_{\q}^ {\prime}(b)} =e^{-\int_{x}^{b} \nu_\q(m)dm}\fr 1{\nu_\q(b)},& x \leq b,\\V^{b}(x)=x-b+ V^{b}(b), & x>b \ec. \la{divL} 
\end{equation}
where the second line  follows upon completing the barrier strategy by "reduce holdings to $b$ when above".

Maximizing over the reflecting barrier $b$ is simply achieved  by finding the roots of
\beq \la{smf} W_\q'' (b^*) =0 \Eq \pd {}{b} \Big[\fr 1{\nu_\q(b^*)} \Big]=1.\eeq
 This is a smooth fit equation at $b^*$ (see \eqr{divL}).\fn[5]{The equivalence between the two   de-Finetti optimality conditions may be checked by differentiating $W_q'(b^*)=W_q(b^*)  \nu_q(b^*)$, which yields $0=W_q''(b^*)=W_q(b^*)(\nu_q'(b^*) + \nu_q^2(b^*))$.}

 Our paper replaces ruin in the  de Finetti dividends barrier optimization by  more general
\AY/\gen \dd  stopping times, to be chosen optimally. Then, the appropriate tool is the calculus of variations/optimal control.
Let us note a  recent related paper using general \dd stopping times  \cite{WangZhou}, who study optimality of barrier policies under a fixed prespecified \dd function.

{\bf First passage theory for spectrally negative Markov processes}.
Prior to \cite{LLZ17b}, the classic and  \dd first passage  literatures   were restricted  mostly to parallel analytic treatments of the two particular cases of  diffusions and of
  \sn \lev processes. \cite{LLZ17b}  showed that a direct
unified approach (inspired by \cite{Leh} in the case of diffusions) may achieve the same results for all time homogeneous Markov processes;  the known results   for diffusions and
  \sn \lev processes are just particular cases of general formulas, once expressed in terms of $W,Z$, or of the differential exit parameters $\nu,\de$
  -- see below.

  Assume  the existence of differential versions of the ruin and \sur \probs:
\begin{Ass}
\label{As}For all $q,\th\geq0$ and $u\leq x$ fixed,  assume  that $\sRui_{\q}^b(x,a)$  and
$\Rui_{q,\th}^b(x,a)$ are differentiable in $b$ at  $b=x$, and in particular that the following limits exist:
\be \la{e:b}
\nu_{q}(x,a)   :=\lim_{\varepsilon\downarrow0}\frac{1-
\sRui_{\q}^{x+\varepsilon}(x,a)}{\varepsilon} \; \bf{(total \; infinitesimal  \; hazard \; rate)}
 \ee
and%
\be \la{e:c}
{\de}_{q,\th}(x,a)   :=\lim_{\varepsilon\downarrow0}
\frac{\Rui_{q,\th}^{x+\varepsilon}(x,a)}{\varepsilon} \;
\bf{(infinitesimal \; spatial\, killing  \; rate)}.
\ee
\end{Ass}
\beR   It turns out that everything reduces
to  the differentiability of the two-sided  ruin and survival probabilities as functions of the upper limit. Informally, we may say that the pillar of first passage  theory for \sn Markov processes is proving the existence of $\nu,\de$.\fn[4]{$\nu$ and $\de$ capture the behavior of excursions  of the process away  from its running maximum.}
  Later, the differential characteristics $\nu$,  $\de$
  were extended in \cite{ALL} to the case of generalized draw-down times, which unify classic  first passage times and \dd times.
\eeR 

Since  results  for \sn \lev  \procs \; (like the \deF \prob \; considered here) require often not much more than the \str , it was natural  to  attempt to extend them
to the \sn strong Markov case. As expected, everything worked out almost smoothly  for ``\lev -type cases" like  random walks \cite{AV},   Markov additive processes \cite{IP}, and L\'evy
processes with $\Omega$ state dependent killing \cite{IP}.\\

    However, diffusions and
  \sn \lev processes were always tackled by different methods until the
 pioneering work \cite{LLZ17b}, who  showed that certain \dd problems could be treated
by a unified approach, inspired by \cite{Leh} in the case of diffusions, which can be extended to all time homogeneous Markov processes.

When switching to \sn \Mar \procs, \;$W_{\q} (x)$ must be replaced by a two variables function
$W_{\q}(x,y)$ (which reduces in the \lev case to $W_{\q} (x,y)=\T W_{\q} (x-y)$, with $\T W_q$ being the scale function of the \lev process).

However,  the existence of $W$, as well as that of the scale function $Z$, are not obvious in the non-\lev case, and
it becomes more convenient  to replace them
  by differential versions $\nu$ and $\de$ defined by\eqr{e:c}, \eqr{e:c} below.
  
  Computing
 $\nu,\de,W,Z$  is still an open problem, even for simple classic processes like the \OU
and the Feller branching diffusion with jumps. However,
    one may cut through this Gordian node by  restricting  to processes for which the limits defining $\nu,\de$  exist,
and leaving to the user the responsibility   to  check this for their process. With this  caveat, the results of \cite{LLZ17b,ALL}
 provide a unifying umbrella for \sn \lev processes, diffusions, branching
\procs (including with immigration), logistic branching
\procs, etc. Surprisingly, all these processes  which were traditionally studied separately, may be viewed as
 particular cases of a unified general \fp theory for \sn \strs !

In this paper we illustrate the potential of this  framework   via one application, a variation of the de Finetti problem of maximizing
expected discounted cumulative dividends, where we replace stopping at ruin by an optimally chosen
Azema-Yor/generalized draw-down stopping time.

{\bf Contents}.  We start by reviewing in Section \ref{s:r}   the classic \sn \fp theory,  and in Section \ref{s:gendd} the \fp theory with
   \gen \dd/\AY stopping times.  Section \ref{s:deFM} introduces the \deF dividends optimization problem with \gen \dd/\AY stopping times for \sn Markov \procs . Section \ref {s:Bol} spells out the
calculus of variations problem to be solved,  and
Section \ref{s:Pon} offers its solution via a Pontryaghin-type approach.

Section 6 presents a detailed analysis of the particular case of \lev \procs. Finally, Section 7 considers a more general class of diffusions (general functions of drifted Brownian motion with particular emphasis on logarithmic cases). 

\sec{Generalized draw-down stopping for processes without positive jumps \la{s:gendd}}
Generalized \dd times appear naturally  in the \AY solution of the Skorokhod embedding problem \cite{AY}, and in the  Dubbins-Shepp-Shiryaev, Peskir  and Hobson optimal stopping problems \cite{dubins1994optimal,peskir1998optimal,hobson2007optimal}.
 Importantly, they allow a unified treatment of classic first passage  and \dd times (see also \cite{ALL} for a further generalization to taxed processes)--see \cite{AVZ,LVZ}.  The idea is to replace  the upper side of the rectangle by a parametrized curve $$(x,y)=({\widehat d}(s),d(s)), \qu {\widehat d}(s)=s-d(s),$$ where $s=x+y$ represents the value of $\ol X_t$ during the excursion which intersects the upper boundary at $(x,y)$. See Figure \ref{f:plDnew}, where  we put ${\widehat d}(s):=s- d(s)$.

 Alternatively, parametrizing by $x$ yields $y=h(x), h(x)=(l)^{-1}(x) -x  $
 (note $Y_t \geq  d(\ol X_t) \Eq  Y_t \geq h(X_t)$).
  \figu{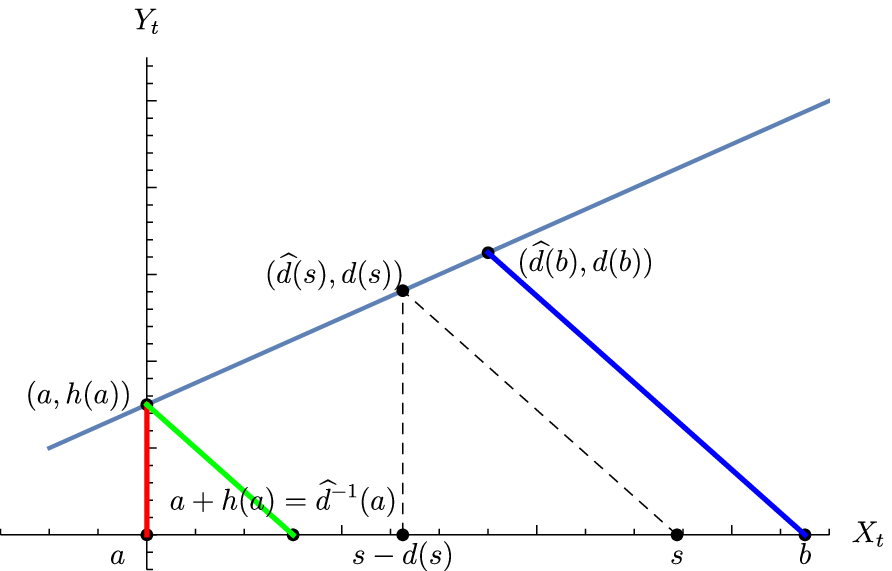}{Exit of $(X,\Y)$ from a trapezoid with $a=0$,
 $d(x)=(1-\xi) x + d=\fr 13 x+ 1$}{0.9}

  \beD \la{d:LVZ}  \cite{AY,LVZ}  For any 
   function $d(s)>0$  such that ${\widehat d}(s)=s-d(s) $ is nondecreasing, a {\bf generalized draw-down}  time is defined by
\be \t_{d}:=\inf\{t\geq0:Y_t>  d(\ol X_t)\}=\inf \left\{t\geq 0: X_t < {\widehat d}(\ol {X}_t)\right\}.
 \label{ddg}\ee
 \iffa
 Introduce $$\T Y_{t}:=Y_t-  d(\ol X_t), \; {t\geq0} $$ 
   to be called  draw-down type
process.
Note that we have $\T Y_{0}=- {\widehat d}(X_{0})<0$,
and that the process $\T Y_t$ is in general  {non-Markovian}.   {However, it is
Markovian during each negative excursion of $X_t$, along one of the oblique lines
in the geometric decomposition sketched in Figure \ref{f:plDnew}.}


\eeD

\beXa  With affine functions
\be d(x)= (1-\xi) x+ d \Eq {\widehat d}(x)= \xi x-d, \xi\in  {[0, 1]} \label{aff}\ee \bea \Eq  h(x)=\fr{(1-\xi) x+ d}{\xi}, \, \xi \in [0,1], d\geq 0, \eea
we obtain the   affine draw-down/regret times studied in \cite{AVZ}.

 Affine draw-down times reduce  to a classic draw-down  time \eqref{dd} when $\xi=1,  d(x)=  d$, and to a ruin time when $\xi=0, {\widehat d}(x)=-d,  d(x)=  x+d$. When  $\xi $ varies, we are dealing with the pencil of lines passing through $(x,y)=(-d,d)$.
In particular, for $\xi=1$ we obtain an infinite strip, and for $\xi=0, d=0,$ we obtain the positive quadrant  (this case corresponds to the classic ruin time).

One of the merits of affine draw-down times is that they allow unifying the classic first passage  theory with the  \dd theory \cite{AVZ}. 
 A second merit is that they {intervene} 
 in the variational problem considered below.

\eeXa

\sec{Optimal dividends barrier problem  for \sn Markov    processes
with \gen draw-down stopping \la{s:deFM}}
Consider now  the extension of de Finetti's optimal dividend problem \be \la{tdivd}
 V^{b}(x)=V_{\q,{\widehat d}(.)}^{b}(x):=  E_x \left[ \int_{0}^{\td^{b]}_d} e^{-q  t} \md(\ol {X}_t-b) \right],\ee
 where $\td^{b]}_d$ denotes a generalized \dd time for the process $X_t^{b]}$ reflected at $b$. Note that $V^{b}$ depends now also on the ``spatial killing function" ${\widehat d}(.)$.
 \beR
 This definition assumes that the initial point satisfies $X_0=\ovl X_0=x$, i.e. that the starting point is on the $x$ axis in figure \ref{f:plDnew}. 
 \eeR

 The  \str \; yields again an explicit decomposition formula
 \be V^{b}(x)=E_{x}\left[  e^{-q\tb
}1_{\left\{  \tb<\t^{b]}_d\right\}  }\right]  V^{b}(b). \la{V} \ee
Furthermore,
by \cite[Thm1]{ALL} 
\ith
 \be
 E_{x}\left[  e^{-q\tb
}1_{\left\{  \tb<\t^{b]}_d\right\}  }\right]=e^{-
\int_{x}^{b}\nu_{q}(z,{\widehat d}(z))\mathrm{d}z}, \la{V1}\ee
where $\nu_{\q}(x,{\widehat d}(x))$ is defined in \eqr{e:b}.

Concerning the expectation of the dividends starting from the barrier $v(b)=V^{b}(b)$, one may show again via standard bounding arguments (see for example \cite[Sec 4]{czarna2018fluctuation}) that 
\be  \la{VL} v(b)=v_{\q}(b,{\widehat d}(b)):=  E_b \left[ \int_{0}^{T_d } e^{-q  t}\md(\ol {X}_t-b) \right]=\nu_{\q}(b,{\widehat d}(b))^{-1}. \ee

Note that in the \lev case, using $x-{\widehat d}(x)=d(x)$, the equations above simplify to:
 \bea
 V^{b}(x)=\fr{W_{\q}(d(x))}{W_{\q}(d(b))} \nu_{\q}(d(b))^{-1},\eea
 which checks with \cite[Lem. 3.1-3.2]{WangZhou}.

\sec{A variational problem for de Finetti's optimal dividends until a generalized \dd time,
with a bound on the initial and total  draw-down/regret area \la{s:Bol}}
Let us consider now de Finetti's optimal dividends
with   draw-down stopping. Suppose $X_0= \ovl X_0=a \geq 0$, and  view  the total trapezoidal area $A(b)=\int_a^b \sqrt{2} d(s) \md s$ between the green and blue lines, in which the bivariate process $(X_t,\Y_t)$ is allowed to evolve, as a measure of risk. With no upper bounds on $A(b)$, the optimum will be $b=\I$.  We set therefore an upper limit
$ A(b) \leq K \sqrt{2}$, and also an upper limit $d_0=d(a)$ on the initial maximum regret.
Using  \eqr{V}-\eqr{VL}, we arrive  to the following Bolza problem
{ \begin{equation} \la{BolV}
\bc \max_{d(y)\geq 0,\
   {d(y), y-d(y) \text{ nondecreasing}}} V ^{b}(a)= \max_{d(y)\geq 0,   {d(y), y-d(y) \text{ nondecreasing}}}
\fr{e^{-\int_a^b \nu_q(y,{\widehat d}(y)) \md y}}{\nu_q(b,{\widehat d}(b))}\\A(b) \leq K \sqrt{2} \\
b,d(b) \text{  free}, d(a) \leq d_0\ec
 \end{equation}


 After taking logarithms, \eqr{BolV}  becomes:
 \begin{equation} \la{Bol} \bc \min_{d(y)\geq 0,   {d(y), y-d(y) \text{ nondecreasing}}} \int_a^b \nu_q(y,{\widehat d}(y)) \md y + \log({\nu_q(b,{\widehat d}(b))}=\\
 \min_{d(y)\geq 0,    {d(y), y-d(y) \text{ nondecreasing}}} \int_a^b \Big(\nu_q(y,{\widehat d}(y))+ \fr{\nu_q'(y,{\widehat d}(y))}{\nu_q(y,{\widehat d}(y))} \Big) \md y + \log({\nu_q(a,{\widehat d}(a))}\\\int_a ^b d(y)  dy
\leq K \\
b,d(b) \text{  free}, d(a) \leq d_0.\ec
 \end{equation}

\beR  Let us relate \eqr{Bol} to the classic \deF problem, which is the particular case obtained by imposing the additional constraint ${\widehat d}(y)={\widehat d}(a) \Eq d(y)=d(a)+y-a \Lra X_t \geq {\widehat d}(a), \for t$. Here the constraint $d(a) = d_0 \Eq {\widehat d}(a) = a -d_0$ quantifies an imposed initial  bankruptcy level, and the subsequent values $d(y)$  quantify   bankruptcy levels dependent on the attained maximum $y$. The area
constraint is thus an acceptable "integrated bankruptcy  risk".
\eeR

\beR If we fix  the \dd boundary in  \eqr{Bol},  the optimality condition for the   dividend barrier $b^*$
is \beq \nu_q(b^*,{\widehat d}(b^*))+ \fr{\nu_q'(b^*,{\widehat d}(b^*))}{\nu_q(b^*,{\widehat d}(b^*))}=0, \la{feq} \eeq
which implies   the classic smooth-fit equation \eqr{smf}.
In the \lev case, the optimal \dd boundary turns out to be De Finetti's ${\widehat d}(y)=const$, and thus the smooth-fit equation at $b^*$ determines  completely the solution.
\eeR

\sec{Solving the \deF Markovian variational problem by Pontryaghin's minimum principle \la{s:Pon}}

    As usual in modern calculus of variations, we let $u(t)$ denote the derivative of $d(t)$, and reformulate the problem
 as
  \begin{equation}
\label{Problem}
\left\lbrace
\begin{split}
&V(a,d(a)):=\inf_{u}J_b(a,d,u),\\
&where\\
&J_b(a,d,u):=\left.\begin{split}&\int_{{a}}^b\pr{\nu_q\pr{t,t-d(t,u)}+
\frac{\partial_1\nu_q\pr{t,t-d(t,u)}+\partial_2\nu_q\pr{t,t-d(t,u)}(1-u(t))}{\nu_q\pr{t,t-d(t,u)}}}dt
\\&{+\log\pr{\nu_q\pr{a,a-d(a)}}}\end{split} \right.\\
&s.t.\ \int_{{a}}^b d(t,u)dt \leq K,\\
&\partial_t d(t,u)=u(t),\ d(a,u)=d(a)\in [0,d_0]\\
& u\textnormal{ measurable }, u \in \pp{u_*=0,u^*=1}.
\end{split}
\right.
\end{equation}
\beR
 In the case of non-decreasing draw-down functions, requiring $d(t)\geq 0$ amounts to imposing the initial condition $d(a)\geq 0$. In absence of such assumptions, one deals with state constraints and Pontryagin's principle has a (slightly) different form
\eeR

This  is the first step towards defining an associated Hamiltonian $H(t,d,u,p)$ \eqr{H1_gen}, where the \textit{costate} $p(t)$ satisfies the {conjugate} equation
\eqr{costate_gen}.  Then, one may apply Pontryagin's maximum principle \cite{pontryagin2018mathematical}.

It is convenient here to break the solution in four cases,  with free area constraint
$$\int_{{a}}^{b^*}d^{opt}(t)dt=(resp.<)K,$$ respectively  with free starting data $d(a)>0$ and with fixed starting data $d(a)=0$.%
\subsection{Optimality Without Area Constraints}
\subsubsection{Arguments for Free Initial  $d(a)>0$}
The associated Hamiltonian {to be minimized} is
\begin{equation}
\label{H1_gen}
\begin{split}
H(t,d,u,p):=&pu+\pr{\nu_q\pr{t,t-d}+\frac{\partial_1\nu_q\pr{t,t-d}+\partial_2\nu_q\pr{t,t-d}(1-u)}{\nu_q\pr{t,t-d}}}\\
=&\nu_q\pr{t,t-d}+\frac{\partial_1\nu_q\pr{t,t-d}+\partial_2\nu_q\pr{t,t-d}}{\nu_q\pr{t,t-d}}+
\pr{p-\frac{\partial_2\nu_q\pr{t,t-d}}{\nu_q\pr{t,t-d}}}u,
\end{split}
\end{equation}
with \textit{costate} $p(t)$ satisfying:
\begin{equation}
\label{costate_gen}
\begin{split}
\partial_tp(t)=&{-\partial_d H(t,d^{opt},u^{opt},p)}\\
=&\partial_2\pr{\nu_q+\frac{\partial_1\nu_q+\partial_2\nu_q}{\nu_q}}\pr{t,t-d^{opt}(t)}-\partial_{2}\Big[\frac{\partial_2\nu_q}{\nu_q}\Big]\pr{t,t-d^{opt}(t)}u^{opt}(t).
\end{split}
\end{equation}
\beR
Note that due to linearity in $u$, optimizing the Hamiltonian $H(t,d,u,p)$ yields optimal control policies $u^{opt}$  of bang-bang type,   except on sets where $p(t)=\frac{\partial_2\nu_q\pr{t,t-d(t)}}{\nu_q\pr{t,t-d(t)}}$.
\eeR
The previous remark implies:
\beL \la{l:3} The optimal draw-down function $d$ may have   three possible types of subintervals.
\BEN \im
On sets $\pp{\alpha_1,\beta_1}$ on which $u^{opt}=u_*=0$, it follows from  {\eqr{Problem}} that $d^{opt}(t)=d(\alpha_1)$ is constant. On such sets, \eqr{costate_gen} yields \begin{equation}
\label{CostateOptCt_gen}
p(t)=p(\alpha_1)+\int_{\alpha_1}^t \partial_2\pr{\nu_q+\frac{\partial_1\nu_q+\partial_2\nu_q}{\nu_q}}\pr{s,s-d\pr{\alpha_1}}ds> \frac{\partial_2\nu_q\pr{t,t-d\pr{\alpha_1}}}{\nu_q\pr{t,t-d\pr{\alpha_1}}},
\end{equation}
for all $t\in\pr{\alpha_1,\beta_1}$ by noting that the coefficient of $u$ in $H$ should be positive.
We either have $\alpha_1=a$ or equality at $t=\alpha_1$ in (\ref{CostateOptCt_gen}). Similar assertions hold true at $t=\beta_1$.
\im  Sets $\pp{\alpha_2,\beta_2}$  on which the costate satisfies the structural equality $p(t)=\frac{\partial_2\nu_q\pr{t,t-d(t)}}{\nu_q\pr{t,t-d(t)}}.$ Recalling that $\partial_td(t)=u(t)$ and combining with (\ref{costate_gen}), {whenever the function $\nu_q$ is regular enough (of class $C^2$ such that second-order mixed partial derivatives coincide)}, this leads to the following implicit ``structural equation" satisfied by the optimal draw-down $d^{opt}$:{
\begin{equation}
\label{Structure_dOpt_gen}
\partial_2\nu_q\pr{t,t-d^{opt}(t)}=0.
\end{equation}}
\im The third and last case to be taken into consideration leads to $u^{opt}=u^*=1$. On sets $\pr{\alpha_3,\beta_3}$ corresponding to this case (note that here \ith that $t-d(t)=\alpha_3-d\pr{\alpha_3}$), one gets \begin{equation}
\label{nonCtCostateOptCt_gen}
\left\lbrace
\begin{split}
d(t)&=d(\alpha_3)+\pr{t-\alpha_3},\\
p(t)&=p\pr{\alpha_3}+\int_{\alpha_3}^t \partial_2\nu_q\pr{s,\alpha_3-d\pr{\alpha_3}}ds+\frac{\partial_2\nu_q}{\nu_q}
\pr{t,\alpha_3-d\pr{\alpha_3}}-{\frac{\partial_2\nu_q}{\nu_q}\pr{\alpha_3,\alpha_3-d\pr{\alpha_3}}}\\
&\leq \frac{\partial_2\nu_q\pr{t,\alpha_3-d\pr{\alpha_3}}}{\nu_q\pr{t,\alpha_3-d\pr{\alpha_3}}}
\end{split}
\right.
\end{equation}
\EEN
\eeL
 \beR
\label{RemarkTransv}
One should add the following transversality conditions, taking into account the liberty to choose $b$ and   $d(b)$, and the initial  conditions $d(a)$:
\BEN \im[i.]   The first condition is linked to the freedom of $b$
\begin{equation}\label{condH_gen}H(b^*)=0.\end{equation}
\item[ii.]  The second condition is linked to the freedom of $d(b)$
\begin{equation}\label{condP_gen}p(b^*)=0.\end{equation}
\item[iii.] If the initial position $d(a)$ is not fixed (thus, one searches for $0<d(a)<d_0$), one further imposes
\begin{equation}\label{condP0_gen}{p(a)=\frac{\partial_2\nu_q}{\nu_q}\pr{a,a-d(a)}}.\\\end{equation}  In other words, assuming optimality, either the optimal initial position satisfies $0<d(a)<d_0$ and, in such cases, (\ref{condP0_gen}) holds true, or, otherwise, $d(a)\in\set{0,d_0}$ (saturating this constraint for some a priori given $d_0$).
\EEN
\eeR

\subsubsection{Arguments for Null Initial Datum $d(a)=0$}
The program presented before still holds true but having fixed the initial datum $d(a)=0$, the transversality conditions are reduced to (\ref{condH_gen}, \ref{condP_gen}).
\subsection{Optimality With Area Constraints}
Again, as before, we reason for $d(a)>0$ (the restriction $d(a)=0$ being taken into account by the absence of the transversality condition  {(\ref{condP0_gen})}.

To cope with the additional constraint $\int_{{a}}^{b^*}d(y)dy{\leq }K$, we use a classic trick and introduce a further variable in the control system. We deal now with
\begin{equation}
\left\lbrace
\begin{split}
&V(a,d(a),e(a)):=\inf_{u}J_b^{{+}}(a,d,e,u),\\
where\ &J_b^{{+}}(a,d,e,u):=J_b(a,d,u),\\
s.t.&\ \partial_t d(t,u)=u(t),\partial_t e(t,u)=d(t,u),\ e(a,u)=e(a),\ d(a,u)=d(a),\ e(b,u)=K.\\
& u\in\pp{0,1}.
\end{split}
\right.
\end{equation}
The associated Hamiltonian is
\begin{equation}
\label{H_2_gen}
\begin{split}
H&^+(t,d,e,u,p,r):=H(t,d,u,p)+rd\\
&=\nu_q\pr{t,t-d}+\frac{\partial_1\nu_q\pr{t,t-d}+\partial_2\nu_q\pr{t,t-d}}{\nu_q\pr{t,t-d}}+\pr{p-\frac{\partial_2\nu_q\pr{t,t-d}}{\nu_q\pr{t,t-d}}}u+rd.
\end{split}
\end{equation}
The arguments are exactly the same but the equations of the \textit{costates} are given here by
\begin{equation}
\label{costate_2_gen}
\left\lbrace
\begin{split}
\partial_t&r(t)=\pr{-\partial_e H^+(t,d,e,u,p,r)}=0 {\Lra r(t)=r=const}\\
\partial_t&p(t)=\pr{-\partial_d H^+(t,d^{opt},e^{opt},u^{opt},p,r)=}-\partial_d H(t,d^{opt},u^{opt},p)-r\\
&=\partial_2\pr{\nu_q+\frac{\partial_1\nu_q+\partial_2\nu_q}{\nu_q}}\pr{t,t-d^{opt}(t)}+
\frac{\pr{\partial_{2}\nu_q}^2-\partial_{22}^2\nu_q\times\nu_q}{\nu_q^2}\pr{t,t-d^{opt}(t)}u^{opt}(t)-r.
\end{split}
\right.
\end{equation}
Cases are exactly the same as before. Formulas (\ref{CostateOptCt_gen}) and (\ref{nonCtCostateOptCt_gen}) are similar:
 \begin{equation}
\label{CostateOptCt_2_gen}
p(t)=p(\alpha_1)+\int_{\alpha_1}^t \partial_2\pr{\nu_q+\frac{\partial_1\nu_q+\partial_2\nu_q}{\nu_q}}\pr{s,s-d\pr{\alpha_1}}ds-r\pr{t-\alpha_1}> \frac{\partial_2\nu_q\pr{t,t-d\pr{\alpha_1}}}{\nu_q\pr{t,t-d\pr{\alpha_1}}},
\end{equation} resp. \begin{equation}
\label{nonCtCostateOptCt_2_gen}
\left\lbrace
\begin{split}
d(t)&=d(\alpha_3)+\pr{t-\alpha_3},\\
p(t)&=p\pr{\alpha_3}+\int_{\alpha_3}^t \partial_2\pr{\nu_q\pr{s,\alpha_3-d\pr{\alpha_3}}}ds+\frac{\partial_2\nu_q}{\nu_q}\pr{t,\alpha_3-
d\pr{\alpha_3}}\\&-{\frac{\partial_2\nu_q}{\nu_q}\pr{\alpha_3,\alpha_3-d\pr{\alpha_3}}-r\pr{t-\alpha_3}}\\
&\leq \frac{\partial_2\nu_q\pr{t,\alpha_3-d\pr{\alpha_3}}}{\nu_q\pr{t,\alpha_3-d\pr{\alpha_3}}}
\end{split}
\right.
\end{equation}
Finally, the  structure equation (\ref{Structure_dOpt_gen}) becomes\fn[4]{Or, in a more symmetric form,
$\fr{W_\q \partial_{12} W_q -\partial_1 W_q \partial_2 W_q}{W_\q^2}\pr{t,t-d^{opt}(t)}=r$.}
\begin{equation}
\label{Structure_dOpt_2_gen}
{\partial_2\nu_q\pr{t,t-d^{opt}(t)}=r.}
\end{equation}
The transversality conditions (see Remark \ref{RemarkTransv}) are similar and allow to determine the  optimal horizon $b^*$.

We have proven
\begin{Thm}
Assume that the problem (\ref{Problem}) admits an optimal pair $(d^{opt},b^*)$ such that $d^{opt}$ is smooth of class $Lip_{u*}(\pp{0,b^*}$, non-decreasing, with $t\mapsto t-d^{opt}(t)$ non-decreasing.
Then,
{\begin{enumerate}
\item The equations \eqr{costate_2_gen}-\eqr{Structure_dOpt_2_gen} provide the three possible behaviors for $d^{opt}$, while $b^*$ is determined from the transversality conditions \eqr{condH_gen}, \eqr{condP_gen} written for the extended Hamiltonian $H^+$ instead of $H$ to which one adds the area saturation $d^{opt}\pr{b^*}=K$;
\item    {In the absence of the area constraint, the equation \eqr{Structure_dOpt_2_gen} holds  with $r=0$.}
\end{enumerate}}
\end{Thm}
\iffa
\beR \begin{itemize}
\item[i.] It can be easily shown that \textit{overoptimizing} in the sense of allowing unbounded derivatives for $d$ by setting $u^*=\infty$ (and not respecting the contraint $y-d(y)$ nondecreasing) leads to rather trivial results: either constant $d$ or $d(\cdot)$ continuously evolving among the points of the critical set for $\nu_q$ (usually void).
\item[ii.] {In the classical L\'evy framework and for a certain class of diffusions, we will give reasonable (and rather general) conditions on $\nu_q$ yielding affine optimal draw-down}.
\end{itemize}

\eeR

\section{Back  to the \lev case: De Finetti's solution is optimal \la{s:lev}}
Let us go back to the \lev case where
\[\nu_q\pr{x,y}=\T \nu_q(x-y).\]
This case has the further particularity that $\pr{\partial_1+\partial_2}\pr{\nu_q}=0$.

In the rest of this section we will drop the tilde in $\T \nu_q(x)$. Recall that the one-variable functions $\nu_q$ is  non-increasing and non-zero.
In this framework, the  (time-homogeneous) extended Hamiltonian (for which we drop the superscript $^+$) is given by \begin{equation}
\label{HL}
H(d,u,p,{r}):=\nu_q\pr{d}+\pr{\frac{\nu_q'(d)}{\nu_q(d)}+p}u{+rd}.
\end{equation}
\Thr like in  any homogeneous setting, the Hamiltonian $H$ is constant  and  by the transversality condition $H(b^*)=0$ it must equal $0$ along optimal trajectories.

The costate (cf. (\ref{costate_2_gen})) \sats
\begin{equation}
\label{costate}
\partial_tp(t)=-\nu_q'\pr{d^{opt}(t)}-\Big[\frac{\nu_q^{'}}{\nu_q}\Big]'\pr{d^{opt}(t)}{u^{opt}(t)}-r
\end{equation}
and the structural equation (\ref{Structure_dOpt_2_gen})  for  {non-extremal solutions} in the present setting writes down
\begin{equation}
\label{Structure_dOpt}
\nu'_q\pr{d^{opt}(t)}={-r}.
\end{equation}
{
We proceed to show now that only sets with $u=1 \Eq d^{opt}(t)=d(a)+t-a \Eq {\widehat d}(t)={\widehat d}(a)$ constant (the \deF solution) are possible in the \lev case.
\begin{enumerate}
\item Sets on which $p(t)=-\frac{\nu_q'}{\nu_q}\pr{d^{opt}(t)}$ cannot exist, and the optimal control is reduced to bang-bang $0/1$, with or without area constraints. Indeed,
     in this case the  Hamiltonian reduces to \[0=H\pr{d^{opt}(t),u^{opt}(t),p(t),r}=\nu_q\pr{d^{opt}(t)}+rd^{opt}(t)\]
     which is impossible since $\nu_q >0$ and $r\geq 0$ by  \eqr{Structure_dOpt}  (since $\nu_q$ is non-increasing).

\item Sets $\pp{\alpha_1,\beta_1}$ on which $u^{opt}=0$ and $d^{opt}(t)=d^{opt}(\alpha_1)$ is constant cannot exist either. Indeed, on such sets, one must have \begin{equation}
\label{CostateOptCt}
p(t)=p(\alpha_1)-(t-\alpha_1)\pr{\nu_q'\pr{d^{opt}(\alpha_1)}+r}\geq {-}\frac{\nu_q'\pr{d^{opt}(\alpha_1)}}{\nu_q\pr{d^{opt}(\alpha_1)}}, \for t \in \pp{\alpha_1,\beta_1},
\end{equation}
Moreover, since the Hamiltonian is null, it follows that \[\nu_q\pr{d^{opt}(\alpha_1)}+rd^{opt}(\alpha_1)=0.\]
\begin{itemize}
\item \underline{Without area constraints.} Here $r=0$, and the previous equality cannot hold. Thus, $0$ control cannot be optimal.
\item \underline{With $0$ initial datum $d^{opt}(\a_1)=0$.} Again, $0$ control cannot be used.

\item \underline{With area constraints and non-zero initial datum.} $r\neq 0$ and it should be picked such that $r=-\frac{\nu_q\pr{d^{opt}(\alpha_1)}}{d^{opt}(\alpha_1)}$.\\
Since $r$ is negative and $\nu_q$ is non-increasing, the inequality in \eqr{CostateOptCt} is strengthened as $t\geq\alpha_1$ increases, hence $u^{opt}=0$ and $d^{opt}(t)=d(\alpha_1)$ is constant for all $t\in\left( \alpha_1,b^*\right]$. But, then, assuming that $b^*>\a_1$, it follows that   $p(b^*)>{-}\frac{\nu_q'\pr{d(\alpha_1)}}{\nu_q\pr{d(\alpha_1)}}\geq 0$ which contradicts the transversality condition \eqr{condP_gen}.
\end{itemize}
\end{enumerate}}

 We have the following more precise result.
\beL
Let $a$ and the initial datum $d(a)$ be given.
 \begin{enumerate}
 \item Then the optimal draw-down (with or without area constraints) is $d^{opt}(t)=d(t) =d(a)+(t-a)$.
\item Without state constraints, $b^*$ should satisfy the de Finetti 'smooth fit' equation \beq \nu_q\pr{d(b^*)}+\frac{\nu_q'\pr{d(b^*)}}{\nu_q\pr{d(b^*)}}=0 \Eq \partial_b \frac{1}{\nu_q \pr{d(b^*)}}=\partial_b V^{b^*}(b)=1.
    \la{deFcr+}\eeq

\item With area constraints, $b^*_{constr}$ is the minimum between $b^*$ and $b^+$, where $b^+\geq a$ is the solution of the equation \begin{equation}\label{NotSaturated_1}
 {\int_a^{b^+} d^{opt}(s) \md s=\int_0^{b^+-a} \Big( d(a) + y \Big) \md y=d(a) (b^+-a)+(b^+-a)^2/2 =K}.\end{equation}

 \im Without area restriction, the best value is obtained for $d(a)$ extreme
 (i.e. $d(a)\in\set{0,d_0}$).
 \end{enumerate}
\eeL
\begin{proof}{
The first assertion on the optimality of $1$-slope $d$ and the last assertion have been provided prior to the Lemma (recall that the optimal control is $u^{opt}=1$). One has
\begin{equation}
\label{nonCtCostateOptCt}
\left\lbrace
\begin{split}
d^{opt}(t)&=d(a)+t-a,\\
p(t)&=p_0-\nu_q\pr{d^{opt}(t)}-\frac{\nu_q'\pr{d^{opt}(t)}}{\nu_q\pr{d^{opt}(t)}}-rt\leq -\frac{\nu_q'\pr{d^{opt}(t)}}{\nu_q\pr{d^{opt}(t)}}.
\end{split}
\right.
\end{equation}
Since the Hamiltonian should be equal to $0$ (see again the transversality condition \eqr{condH_gen}, it follows that \[p(t)=-\pr{\nu_q\pr{d^{opt}(t)}+\frac{\nu_q'}{\nu_q}\pr{d^{opt}(t)}+rd^{opt}(t)}.\]
We focus on the case without area constraints i.e. $r=0$. The transversality condition \eqr{condP_gen} yields $p\pr{b^*}=0$ which, given the previous form for $p$ yields the second assertion.\\
For the third assertion, we note that the presence of a further constraint (on the area) can only increase the value function. This area is given exactly by $\int_a^{b^*_{constr}} d^{opt}(s) \md s\leq K$. If $b^*\leq b^+$, then it satisfies the constraint, thus providing the best solution. Otherwise, one retains $b^+$.}
\item[(b)] Without area restrictions, the optimal initial datum $d(a)$ is either $d_0$ or $0$ since  if $d(a)$ does not satisfy this restriction, then by the transversality condition \eqr{condP0_gen} it follows that $p(a)=-\frac{\nu_q'}{\nu_q}(d(a))$ and, thus $H\pr{d^{opt}(a),u^{opt}(a),p(a),0}=\nu_q\pr{d^{opt}(a)}>0$. This contradicts the previous assertion on $H$ being $0$.
\end{proof}
\beR
\label{RemarkTransv}
If, instead of searching for draw-down functions s.t. $y\mapsto y-d(y)$ is non-decreasing (i.e. $u^*=1$) one searches for $d^{opt}$ with $u^*<1$-bounded derivative, then the condition \eqr{deFcr+} becomes  \beq \nu_q\pr{d(b^*)}+\frac{\nu_q'\pr{d(b^*)}}{\nu_q\pr{d(b^*)}}u^*=0 \Lra \frac{\nu_q'\pr{d(b^*)}}{\nu_q^2\pr{d(b^*)}}=-   {\fr 1 {u^*} \Eq \partial_b \frac{1}{\nu_q \pr{d(b^*)}}=\fr 1 {u^*}}. \la{deFcru}\eeq

The solution in this case will still be to use $u=u^*$, leading   to the  affine \dd barriers already
 studied in \cite{AVZ} under the different parametrization $d(y)=(1-\xi) \pr{y {-a}}+ d \Eq u^*=1-\xi\in [0,1]$.
 
\eeR

\beXa [\cite{AVZ}] For Brownian motion
with drift  $X_t = \s B_t + \mu  t$, the scale function  is
\be \la{sfB}
W_\q(x) = \frac{1}{\Delta} [
e^{(-\mu+\Delta)x/\s^2}-e^{-(\mu+\Delta) x/\s^2}]:=\frac{1}{\Delta} [
e^{\Phi_\q x}-e^{-\r_\q x}],
\ee
where $\Delta = \sqrt{\mu^2 + 2\q\s^2}$, and  $\Phi_\q,-\r_\q$ are the nonnegative and negative roots of $\s^2 \th^2 + 2 \mu \th -2 \q=0$.  

 In the case of affine optimal profiles and  with the extra restriction $u^*:=1-\xi\leq 1$ as in Remark \ref{RemarkTransv}, (ii) the transversality conditions \eqr{deFcru} yield

{\be \la{xiopt} \frac{\pr{W_{q}''W_{q}-\pr{W_{q}'}^2}}{\pr{W_{q}'}^2}(d(b^*))=-\fr 1 {1-\xi} \Lra \fr{W_{\q}'' W_{\q}}{(W_{\q}')^2}(d(b^*))=- \fr{\xi}{1-\xi},\ee}
a result already obtained in \cite{AVZ}.\fn[4]{When $\xi=d=0$, we recover in the \cP case the equation $W_{\q}''(b)=0$.}
\eeXa  

\sec{Optimal dividends for functions of a \lev process \la{s:gbm}}

Consider a process implicitly defined by \beq \la{imp} F(X_t)-F(x_0)=Z_t, Log \Big(E_0[e^{s Z_t}]\Big)=t \k(s)\eeq  for an arbitrary increasing \fun \ $F(x)$ and $x_0=0$ (w.l.o.g.). This class of processes generalizes the geometric Brownian motion, obtained when  $F(x)=\ln(\fr{\a}{\b} x + 1), Z_t=\s B_t+ \mu t, \s >0$.

\ssec{Optimal dividends for functions of Brownian motion with drift \la{s:gbm}}
 The
monotone harmonic functions are $\f_\pm(x)=\f^Z_\pm(F(x))$, where $\f^Z_\pm$ are the
monotone harmonic functions of $Z_t$.  The $q$-scale and excursions function may \thr be expressed
in terms of the corresponding one-dimensional characteristics of $Z_t$: \bea && W_q(x,y)=W_q^Z(F(x)-F(y)):=\omega_q(F(x)-F(y)), \\&& \nu_q(x,y)=F'(x) \nu^Z_q(F(x)-F(y)):=F'(x) \mu_q(F(x)-F(y)).\eea

As \wk, \; $\f^Z_\pm(F(x))=e^{r_\pm x}$, where $r_\pm$ are the positive /negative roots  of $\fr{\s^2 }2 r^2 + r \mu -\q =0,$
and   the $q$-scale function is:
$$\omega_q(x):=e^{r_+x}-e^{r_-x}.$$

Recalling that   $\omega_q(x)$ satisfies the equation
\be \fr{\s^2}2 \omega_q^{''}+\mu \omega_q'- q \omega_q =0\Eq \omega_q^{''}= {-}\fr{2 \mu}{\s^2}\omega_q'+\frac{2\q}{\s^2}
\omega_q\ee
\wf \beq \la{mueq} \mu_q'=\fr{w_q ''}{w_q} -\mu_q^2=  {-}\frac{2 \mu}{\s^2}\mu_q+\frac{2\q}{\s^2} -\mu_q^2.\eeq

 \beL \la{l:ggbm} Let $X_t$ be defined by \eqr{imp}, with $F(x)$ strictly \inc,  and set $\D=F( x)-F(y)$. Then:

 A) The structure equation \eqr{Structure_dOpt_2_gen} becomes
  \begin{equation}  \label{StructGBMG} {\mu_q^2\pr{\D^{opt}(x)}} {+}\frac{2 \mu}{\s^2}\mu_q \D^{opt}(x)-\frac{2q}{\s^2}=
  \frac{r}{F'(x)F'\pr{F^{-1}\pr{F(x)-\D^{opt}(x)}}},
 \end{equation}for some $r\geq 0$.\\

B) If $F(x)$ is \frt convex, then  for fixed $x$ the  equation \eqref{StructGBMG}
admits   {exactly} one  solution $\D^{opt}(x)$.
  \eeL

{\bf Proof:}
A) The reader is invited to note that
\begin{equation}
\partial_2\nu_q(x,y)=-
F'(x)F'(y)\mu_q'\pr{\D}=-F'(x)F'\pr{F^{-1}\pr{F(x)-\D}}\mu_q'\pr{\D}\\
\end{equation}

The structure equation becomes \thr
\begin{equation}
\label{StrGBMG}
r=
-F'(y^{opt}(x))F'(x)\mu_q'\pr{\D^{opt}(x)}=-F'(x)F'\pr{F^{-1}\pr{F(x)-\D^{opt}(x)}}
\mu_q'\pr{\D^{opt}(x)}
\end{equation}
(recall the case $r=0$ corresponds to the absence of area constraints).
Dividing  now (\ref{StrGBMG}) by $F'(y)F'(x)$  and substituting \eqr{mueq} yields
\eqref{StructGBMG} leads to our first assertion.\\

For notation simplicity we drop the dependence $\D^{opt}(x)$ and  write $\D$ from now on. \\

B) For the uniqueness assertion, fix $x$ and  assume that $\D$ satisfies \eqref{StructGBMG}. Note now that the applications
\[\D\mapsto \Bigg(\mu_q^2(\D) {+}\fr{2 \mu}{\s^2}\mu_q(\D)-\fr{2 \q}{\s^2}\Bigg) \text{ and }\D\mapsto \pr{\frac{r}{F'(x)F'\pr{F^{-1}\pr{F(x)-\D}}}}\] are {decreasing} with range $[0, \I)$ and increasing
with range $[\frac{r}{\pr{F'(x)}^2}, \frac{r}{{F'(x)} F'(0)}]$, respectively.

Indeed, the derivative of the first
\bea 2 \mu_q'(\D)\Big(\mu_q(\D) {+}\fr{ \mu}{\s^2}\Big)\eea
is negative by the strict monotony of $\mu_q$. The values start from $\I$ since the positivity of $\s$ implies $\nu(0)=\I$, and their positivity follows from the \wk \; \eqr{nub}.

For the second term, note that besides the  negative sign and the inversion, it consists of a  composition of the increasing functions $F'$ (here  convexity of $F$ is used), $F^{-1}$
and $F$\fn[4]{Equivalently, the function $\frac{r}{{F'(x)} F'(y)}$ is decreasing in $y$}. This terms is thus increasing and the   assertion follows.  \qed \\

\ssec{Geometric (logarithmic) Brownian motion \la{s:gbm}}
Consider the diffusion defined by  the SDE \beq \la{gBMS} \fr{d X_t}{\a X_t + \b} = dt + \varepsilon d B_t,\eeq
with coefficients
$$\mu(x)=\a x+ \b, \s(x)=\e(\a x+ \b), \b > 0.$$
By Ito's formula, this process may be represented as:
\beq \la{gBM} &&X_t=(x_0 + \fr {\b} {\a}) e^{Z_t}- \fr {\b} {\a},\;  Z_t= \varepsilon \a B_t+ (\a -\fr{\varepsilon^2 \a^2}2)t \Eq\\&& Z_t=F(X_t), \; F(x)=\ln(\fr{\a x+ \b}{\a x_0+ \b}). \no\eeq

\beR Note that $F$ is  concave and \thr  Lemma \ref{l:ggbm} may  not apply. Numeric experiments reveal
however that a unique solution $\D(x)$ exists sometimes.
\eeR

 The
monotone harmonic functions are $\f_\pm(x)=(\a x+ \b)^{r_\pm},$  where
$$r_\pm=
\frac{1-\frac 2{\varepsilon^2\a}\pm\sqrt{\pr{1-\frac{2}{\varepsilon^2\a}}^2+\frac{8q}{\varepsilon^2\a^2}}}{2}$$
 are the positive /negative roots  of $\fr{\e^2 \a^2}2 r^2 + r(\a -\fr{\e^2 \a^2}2) -\q =0$,
appearing in the  scale function $\omega_q(x):=e^{r_+x}-e^{r_-x}$ associated to the drifted Brownian motion $Z_t= \varepsilon \a B_t+ (\a -\fr{\varepsilon^2 \a^2}2)t$.

The two variables scale function  satisfying $W_q(x,y)=0$
is
$$W_q(x,y)=(\fr{\a x+ \b}{\a y+ \b})^{r_+}-(\fr{\a x+ \b}{\a y+ \b})^{r_-}=\omega_q\pr{\ln\frac{\alpha x+\beta}{\alpha y+\beta}},$$
and its logarithmic derivative is
 \[\nu_q(x,y)=\frac{\partial_1 W_q(x,y)}{W_q(x,y)}=\frac{\alpha}{\alpha x+\beta}\mu_q\pr{\ln\frac{\alpha x+\beta}{\alpha y+\beta}},\]
where $ \mu_q:=\frac{\omega_q'}{\omega_q}.$
Since \[
F'(x)F'(y)= F'(x)F'\pr{F^{-1}\pr{F(x)-\D}}=\frac{\a}{\a x+\b}\frac{\a}{e^{F(x)-\D}}=\pr{\frac{\a}{\a x+\b}}^2e^{\D},\]
Lemma \ref{l:ggbm} A) yields here  \beq  {\mu_q^2(\D)+\pr{\frac{2}{\varepsilon^2\a}-1}\mu_q(\D)-\frac{2q}{\varepsilon^2\a^2}=
r\pr{x+\frac{\b}{\a}}^2e^{-\D} \Lra} \eeq

 \[\frac{r_+e^{r_+\D}-r_-e^{r_-\D}}{e^{r_+\D}-e^{r_-\D}}=
\frac{1-\frac{2}{\varepsilon^2\a}+\pp{\pr{\pr{1-\frac{2}{\varepsilon^2\a}}^2+
\frac{8q}{\varepsilon^2\a^2}+4r\pr{x+\frac{\b}{\a}}^2e^{-\D}}}^{\frac{1}{2}}}{2},\]}
or, equivalently,
\begin{equation}\label{StrucExp_gbm}{r_+e^{r_+\D}-r_-e^{r_-\D}=
\frac{1-\frac{2}{\varepsilon^2\a}+\pp{\pr{\pr{1-\frac{2}{\varepsilon^2\a}}^2+\frac{8q}
{\varepsilon^2\a^2}+4r\pr{x+\frac{\b}{\a}}^2e^{-\D}}}^{\frac{1}{2}}}{2}\pr{e^{r_+\D}-e^{r_-\D}}.}\end{equation}

Finally
\be  d^{opt}(x)=\pr{x{+}\frac{\beta}{\alpha}}\pr{1-e^{-\D(x)}}.\ee
\beR
{Note that if $r=0$ \eqr{StrucExp_gbm} becomes $r_+e^{r_+\D}-r_-e^{r_-\D}=
r_+\pr{e^{r_+\D}-e^{r_-\D}} \Lra r_+=r_-$, which is impossible. With  $r>0$ how, if
an adequate solution $\D$ exists, it is a non-constant function of the position $x$.}
\eeR

In conclusion
\beL

{
\begin{enumerate} Assuming $u \in [0,1-\xi]$, \ith
\item Without area restrictions for geometric Brownian motions, the structure equation (for $r=0$) has no solution and the optimal profile is still affine as before with the maximal slope $1-\xi$.
\item For area restrictions, the structure equation is given by \eqr{StrucExp_gbm}. The optimal draw-down belongs to the class of  functions whose piecewise components either satisfy \eqr{StrucExp_gbm} (for some fixed $r>0$) or are affine with the maximal slope $1-\xi$.
\end{enumerate}}
\eeL

{{\it Acknowledgement}}.  We thank Hongzhong Zhang for help in formulating the variational problem.

\small
\bibliographystyle{alpha}

\bibliography{Pare37}
\end{document}